\documentclass[3p,times,procedia]{elsarticle}
\usepackage{ecrcR,fleqn} 
\volume{17} 
\firstpage{1134}
\journalname{Procedia Computer Science}
\runauth{Pavel Chebotarev and Rafig Agaev}  
\jid{procs} 
\jnltitlelogo{Procedia Computer Science}  

\usepackage{amssymb} 
\usepackage{amsthm}  

\usepackage[figuresright]{rotating} 

\def\widetilde{\tilde}
\newcommand\condition[2] {\medskip\noindent{\bfseries #1.\ } {\rm #2}\par\medskip}  
\newcommand\conditiont[2]{\noindent{\bfseries #1.\ }         {\rm #2}\par\medskip}  
\def\e{w}                                                
\def\G{\Gamma}                                           
\def\si{\sigma}                                          
\def\1n{1,\ldots,n}                                      
\def\0n{0,\ldots,n}                                      
\def\GG{\mathop{\mathcal G}\nolimits}                    
\def\FF{\mathop{\mathcal F}\nolimits}                    
\def\TT{\mathop{\mathcal T}\nolimits}                    
\def\PP{\mathop{\mathcal P}\nolimits}                    
\def\interca{\mathop{\scriptscriptstyle{\rm T}}\nolimits}
\def\cdc{,\ldots,}                                       
\def\ktil {\widetilde K}                                 
\def\_#1{\mathop{\hspace{-2pt}^{}_{#1}}}                 
\def\vj {\mathop{\widetilde{J}}\nolimits}                
\def\q  {\widetilde{J}}                                  
\def\J  {\widetilde{J}{}}                                
\def\x{{}}                                               
\def\beq{\begin{equation}}                               
\def\eeq{\end{equation}}                                 
\def\suml {\mathop{\sum}   \limits}                      
\def\l{\ell}                                             
\newcommand{\card}[1]{\left|#1\right|}                   
\def\adj{\mathop{\rm adj}\nolimits}                      
\def\tor{\to\hspace{-.03em}{*}\hspace{-.05em}}            
\def\rto{\hspace{-.00em}{*}\hspace{-.02em}\to}           
\def\di{d'}                                              
\def\ve{v}                                               
\def\rank{\mathop{\rm rank}\nolimits}                    
\def\a{\mathop{\alpha}\nolimits}                         
\def\id{\mathop{\rm id}\nolimits}                        
\def\od{\mathop{\rm od}\nolimits}                        
\def\tr{\mathop{\rm tr}\nolimits}                        
\def\D{{\rm\Delta}}                                      
\def\xy{\hspace{.07em}}

\newtheorem{theorem}{Theorem}
\newtheorem{corollary}{Corollary}

\newtheorem{prop}{Proposition}

\newtheorem{definition}{Definition}{\bfseries}{\bfseries}

\begin{document}
\begin{frontmatter}



\dochead{Information Technology and Quantitative Management , ITQM 2013}

\title{Matrices of forests, analysis of networks, and ranking problems}

\author[a]{Pavel Chebotarev\corref{cor1}}
\author[a]{Rafig Agaev}

\address[a]{Institute of Control Sciences of the Russian Academy of Sciences, 65 Profsoyuznaya str., Moscow 117997, Russia}

\begin{abstract}
The matrices of spanning rooted forests are studied as a tool for analysing the structure of networks and measuring their properties. The problems of revealing the basic bicomponents, measuring vertex proximity, and ranking from preference relations / sports competitions are considered. It is shown that the vertex accessibility measure based on spanning forests has a number of desirable properties. An interpretation for the stochastic matrix of out-forests in terms of information dissemination is given.
\vskip11pt
\end{abstract}

 \cortext[cor1]{Corresponding author. Tel.: +7-495-334-8869; fax: +7-495-420-2016.}
\end{frontmatter}
\email{upi@ipu.ru; pavel4e@gmail.com}


\section{Introduction}

The matrices of walks between vertices are useful to analyse the structure of networks (see, e.g., \cite{Che12DAM} and the references therein). These matrices are the powers of the adjacency matrix. In this paper, we consider the matrices of spanning rooted forests as an alternative tool for analyzing networks (cf.~\cite{SenelleGarcia-DiezMantrachShimboSaerensFouss13X}). We show how they can be used
for measuring vertex proximity (Section~\ref{dosti}) and for ranking on the base of preference relations / sports competitions (Section~\ref{leader}). In the first sections of the paper, we introduce the necessary notation (Section~\ref{Notatio}) 
and list some properties of spanning rooted forests and forest matrices (Section~\ref{sect_prop}).

Three features that distinguish the matrices of forests from the matrices of walks are notable. First, all column sums (or row sums) of the forest matrices are the same, therefore, these matrices can be considered as matrices of {\em relative\/} accessibility. Second, there are matrices of ``out-forests'' and matrices of ``in-forests'', enabling one to distinguish ``out-accessibility'' from ``in-accessibility'', which is intuitively justifiable. Third, the total weights of maximum spanning forests are closely related to the Ces\'aro limiting probabilities of Markov chains determined by the network under consideration.

\section{Notation and simple facts}
\label{Notatio}

\subsection{Networks, components, and bases}
\label{sec2}

Suppose that $\G$ is a weighted digraph (= network) without loops, $V(\G)=\{\1n\},$ $n>1,$ is its set of vertices and $E(\G)$ the set of arcs. Let $W=(\e\_{ij})$ be the matrix of arc weights. Its entry $\e\_{ij}$ is zero if there is no arc from vertex $i$ to vertex~$j$ in~$\G$; otherwise $\e\_{ij}$ is strictly positive.

The {\it column Laplacian\/} matrix of $\G$ is the $n\times n$ matrix $L=L(\G)=(\l\_{ij})$ with entries
$\l\_{ij}=-\e\_{ij}$ whenever $j\ne i$ and $\l\_{ii}=-\sum_{k\ne i}\l\_{ki}$, $i,j=\1n$. This matrix has
zero column sums. 
The column Laplacian matrices are singular M-matrices (see, e.g., \cite[p.~258]{CheAga02a+}). Their index is 1~\cite[Proposition~12]{CheAga02a+}.

If $\G'$ is a subgraph of $\G$, then the weight of $\G'$, $\e(\G')$, is the product of the weights of all its arcs; if $E(\G')=\varnothing$, then $\e(\G')=1$ by definition. The weight of a nonempty set of digraphs $\GG$ is
\beq
\label{set_weight}
\e(\GG)=\suml_{H\in\GG}\e(H);\quad \e(\varnothing)=0.
\eeq

{\it A spanning subgraph\/} of $\G$ is a subgraph with vertex set $V(\G)$. The {\it indegree\/} $\id(\ve)$ and {\it outdegree\/} $\od(\ve)$ of a vertex $\ve$ are the number of arcs that come {\em in} $\ve$ and {\em out of} $\ve$, respectively. A vertex $\ve$ is called a {\it source\/} if $\id(\ve)=0$. A vertex $\ve$ is {\it isolated\/} if $\id(\ve)=\od(\ve)=0$. A {\it walk\/} ({\it semiwalk}) is an alternating sequence of vertices and arcs $\ve\_0,e\_1,$ $\ve\_1\cdc e\_k, \ve\_k$ with every arc $e\_i$ being $(\ve\_{i-1},\ve\_i)$ (resp., either $(\ve\_{i-1},\ve\_i)$ or $(\ve\_i,\ve\_{i-1})$). A {\it path\/} is a walk with distinct vertices. A~{\it circuit\/} is a walk with $\ve\_0=\ve\_k$, the other vertices being distinct and different from $\ve\_0$. Vertex $\ve$ {\it is reachable\/} from vertex $z$ in $\G$ if $\ve=z$ or $\G$ contains a path from $z$ to~$\ve$.

A digraph is {\it strongly connected\/} (or {\it strong}) if all of its vertices are mutually reachable and {\it weakly connected\/} if any two different vertices are connected by a semiwalk. Any maximal strongly connected (weakly connected) subgraph of $\G$ is a {\it strong component}, or a {\it bicomponent} (resp., a {\it weak component}) of~$\G$. Let $\G_1\cdc\G_r$ be all the strong components of~$\G$. The {\it condensation\/} (or {\it factorgraph}, or {\it leaf composition}, or {\it Hertz graph}) $\G^{\circ}$ of digraph $\G$ is the digraph with vertex set $\{\G_1\cdc\G_r\}$, where arc $(\G_i,\G_j)$ belongs to $E(\G^{\circ})$ iff $E(\G)$ contains at least one arc from a vertex of $\G_i$ to a vertex of~$\G_j$. The condensation of any digraph $\G$ obviously contains no circuits.

A {\em vertex basis\/} of a digraph $\G$ is any minimal (by inclusion) collection of vertices such that every vertex of $\G$ is reachable from at least one vertex of the collection. If a digraph does not contain circuits, then its vertex basis is obviously unique and coincides with the set of all sources~\cite{Harary69,Zykov69}. That is why the bicomponents of $\G$ that correspond to the sources of $\G^{\circ}$ are called the {\it basic bicomponents\/}~\cite{Zykov69} or {\it source bicomponents\/} of~$\G$. In this paper, the term {\it source knot of\/}~$\G$ will stand for the set of vertices of any source bicomponent of~$\G.$ In~\cite{FiedlerSedlacek58}, source knots are called {\em W-bases}.

The following statement \cite{Harary69,Zykov69} characterizes all the vertex bases of a digraph.

\begin{prop}
\label{proZy} A set $U\subseteq V(\G)$ is a vertex basis of\/ $\G$ if and only if\/ $U$ contains exactly one vertex from each source knot of\/ $\G$ and no other vertices.
\end{prop}

Schwartz \cite{Schwartz86} referred to the source knots of a digraph as the {\em minimum $P$-undominated sets}. According to his {\em Generalized Optimal Choice Axiom\/} (GOCHA), if a digraph represents a preference relation on a set of alternatives, then the {\it choice\/} should be the union of its minimum $P$-undominated sets.\footnote{This union is also called the {\em top cycle\/} and the {\em strong basis\/} of~$\G$.} This choice is interpreted as the set of ``best'' alternatives. A review of choice rules of this kind can be found in \cite{Vol'skii88}; for ``fuzzy'' extensions, see~\cite{Roubens96FSS}.

\subsection{Matrices of forests}

A {\it diverging tree\/} is a weakly connected digraph in which one vertex (called the {\it root}) has indegree zero and the remaining vertices have indegree one. A~diverging tree is said to {\em diverge from\/} its root. Spanning diverging trees are sometimes called {\it out-arborescences}. A~{\it diverging forest\/} (or {\it diverging branching}) is a digraph all of whose weak components are diverging trees. The roots of these trees are called the roots of the diverging forest. A~{\em converging tree\/} ({\em converging forest}) is a digraph that can be obtained from a diverging tree (resp., diverging forest) by the reversal of all arcs. The roots of a converging forest are its vertices that have outdegree zero. In what follows, spanning diverging forests in $\G$ will be called {\it out-forests\/} of $\G$; spanning converging forests in $\G$ will be called {\it in-forests\/} of~$\G$.

\begin{definition}\label{De-Max}
{\em An out-forest $F$ of a digraph $\G$ is called a {\em maximum out-forest\/} of $\G$ if $\G$ has no out-forest with a greater number of arcs than in~$F$.}
\end{definition}

It is easily seen that every maximum out-forest of $\G$ has the minimum possible number of diverging trees; this number will be called the {\it out-forest dimension\/} of $\G$ and denoted by~$\di$. It can be easily shown that the number of arcs in any maximum out-forest is $n-\di$; 
the number of weak components in a forest with $k$ arcs is $n-k$.

By $\FF^{\rto}(\G)=\FF^{\rto}$ and $\FF^{\rto}_k(\G)=\FF^{\rto}_k$ we denote the set of all out-forests of $\G$ and the set of all out-forests\x{} with $k$ arcs of $\G$, respectively; $\FF^{i\rto j}_k$ will designate the set of all out-forests with $k$ arcs where $j$ belongs to a tree diverging from~$i$; $\FF^{i\rto j}=\bigcup_{k=0}^{n-\di}\FF^{i\rto j}_k$ is the set of such out-forests with all possible numbers of arcs. The notation like $\FF^{\rto}_{(k)}$ will be used for sets of out-forests that consist of $k$ trees, so $\FF^{\rto}_{(k)}=\FF^{\rto}_{n-k},\;k=\1n.$ Thus, the ${*}\hspace{-.4em}\to$ sign relates to out-forests; the corresponding notation with $\to\hspace{-.3em}{*}\,$, such as $\FF^{\tor},$ relates to in-forests, i.e., $\ast$ images the root(s).
Let
\beq\label{sik}
\si\_k
=\e(\FF^{\rto}_k),\quad k=0,1,\ldots;\quad\si
=\e(\FF^{\rto})=\suml_{k=0}^{n-\di}\si\_k.
\label{si}
\eeq

By (\ref{sik}) and (\ref{set_weight}), $\si\_k=0$ whenever $k>n-\di;$ $\si\_0=1.$
We also introduce the parametric value
\beq\label{sitau}
\si(\tau)
=\suml_{k=0}^{n-\di}\e(\FF^{\rto}_k)\,\tau^k
=\suml_{k=0}^{n-\di}\si\_k\tau^k,\quad \tau>0,
\eeq
which is the total weight of out-forests in $\G$ provided that all arc weights are multiplied by~$\tau$.

Consider the {\em matrices $Q\_k=(q_{ij}^k),\; k=0,1,\ldots,$ of out-forests\x{} with $k$ arcs of\/ $\G$}: the entries of $Q\_k$ are
\beq\label{qijk}
q_{ij}^k=\e(\FF_k^{i\rto j}).
\eeq

By (\ref{qijk}) and (\ref{set_weight}), $Q\_k=0$ whenever $k>n-\di;$
$Q\_0=I.$

The {\em matrix of all out-forests\/} is
\beq\label{qij}
Q=(q\_{ij})=\suml_{k=0}^{n-\di}Q\_k\mbox{~~with entries~~}q\_{ij}=\e(\FF^{i\rto j}).
\eeq

We will also consider the {\em stochastic matrices of out-forests}:
\beq\label{Jk}\label{J}
J\_k=\si_k^{-1}Q\_k,\quad k=\0n-\di;\quad
J=(J\_{ij})=\si^{-1}Q
\eeq
and the parametric matrices
\beq\label{Qtau}\label{Jtau}
Q(\tau)=\suml_{k=0}^{n-\di} Q\_k\tau^k\mbox{  and  }
J(\tau)=\si^{-1}(\tau)\,Q(\tau),\quad \tau>0,
\eeq
where $\si\_k,$ $\si,$ and $\si(\tau)$ are defined by (\ref{sik}) and (\ref{sitau}). $J\_k$ are column stochastic by Theorem~\ref{th4} below.

The {\em stochastic matrix of maximum out-forests\/} $J\_{n-\di}\/$ will also be denoted by~$\J=(\J\_{ij})$:
\beq\label{Jbar}
\J=J\_{n-\di}.
\eeq

The matrices of forests can be found by means of matrix analysis~\cite[Section 5]{CheAga02a+}.

\section{Properties of the forest matrices}
\label{sect_prop}

A number of results on the forest matrices are presented in~\cite{CheAga02a+}. Some of them are collected in the following theorem.

\begin{theorem}\label{th4}\label{sumsum7}\label{teo.allk}
$\!\!\!${\rm~\cite{CheAga02a+,CheSha97}.}
{\rm 1.} Matrices $J\_k,\;k=\0n-\di,$ $J,$ and $J(\tau)$ are column
         stochastic.\\
{\rm 2.} The matrix-forest theorem: for any $\tau>0,$ $Q(\tau)=\adj(I+\tau L)$ and
         $\si(\tau)=\det(I+\tau L),$ 
         whence$,$
         $J(\tau)=(I+\tau L)^{-1}.$\\
{\rm 3.} $L\q=\q L=0$.\\
{\rm 4.} $\vj$ is idempotent$:$ $\;\q^{2}=\q.$\\
{\rm 5.} $\J
         =\lim_{\tau\to\infty}J(\tau)
         =\lim_{\tau\to\infty} (I+\tau\,L)^{-1}.$\\
{\rm 6.} $\rank\J=\di;\;\rank L=n-\di$.\\
{\rm 7.} ${Q\_{k}
         =\sum_{i=0}^k\si\_{k-i}(-L)^i,\;\;\; k=0,1,\ldots.}$\\
{\rm 8.} $\J$ is the eigenprojection of $L$.
\end{theorem}


To formulate the topological properties of the matrix $\J$, the following notation is needed.

Let $\ktil=\bigcup_{i=1}^{\di} K\_i$, where $K\_i$ are all the source knots of~$\G$; let $K_i^{+}$ be the set of all vertices reachable from $K\_i$ and unreachable from the other source knots. For any $k\in\ktil$, $K(k)$ will designate the source knot that contains~$k$. For any source knot $K$ of $\G,$ denote by $\G_K$ the restriction of $\G$ to $K$ and by $\G_{-K}$ the subgraph with vertex set $V(\G)$ and arc set $E(\G)\setminus E(\G_K)$. For a fixed $K$, $\TT$ will designate the set of all spanning diverging trees of $\G_K$, and $\PP$ the set of all maximum out-forests of $\G_{-K}$. By $\TT^k,$ $k\in K,$ we denote the subset of $\TT$ consisting of all trees that diverge from $k$, and by $\PP^{K\rto j},$ $j\in V(\G),$ the set of all maximum out-forests of $\G_{-K}$ such that $j$ is reachable from some vertex that belongs to~$K$ in these forests. $\q_{k\bullet}$ is the $k$th row of~$\q$.

\begin{theorem} 
\label{th2}
$\!\!\!${\rm~\cite{AgaChe00}.}
{Let $K$ be a source knot in~$\G$. Then the following statements hold.\\
{\rm 1.}~$\q_{ij}\ne 0\;\Leftrightarrow\; (i\in\ktil$
         and $j$ is reachable from $i$ in~$\G).$\\
{\rm 2.}~Let $k\!\in\!K.\!$ For any $j\in V(\G),$
         $\q_{kj}=\e(\TT^k)\xy\e(\PP^{K\rto j})\slash
         \e(\FF^{\rto}_{(\di)}).\!$ Furthermore$,$ if $j\in K^+,$ then
         $\q_{kj}=\q_{kk}=\e(\TT^k)\slash \e(\TT).$\\
{\rm 3.}~$\suml_{k\in K}\q_{kk}=1.$ In particular$,$
         if $k$ is a source$,$ then $\q_{kk}=1.$\\
{\rm 4.}~For any $k\_1,k\_2\in K$, $\q_{k\_2\!\bullet}=(\e(\TT^{k\_2})
         \slash\e(\TT^{k\_1}))\xy\q_{k\_1\!\bullet}$ holds$,$ i.e.$,$ the
         rows $k\_1$ and $k\_2$ of $\q$ are proportional.}
\end{theorem}

Thus, knowing $\q$ enables one to reveal basic bicomponents and the digraph reachability structure (cf.~\cite{AgaChe00,Cheb0508:Matrices}).

We say that a weighted digraph $\G$ and a finite homogeneous Markov chain with transition probability matrix $P$ {\em inversely correspond\/} to each other if \beq \label{G_M_cor} I-P=\a L^{\interca}, \eeq where $\a$ is any nonzero real number.

If a Markov chain inversely corresponds to $\G,$ then the probability of transition from $j$ to $i\ne j$ is proportional to the weight of arc $(i,j)$ in $\G$ and is $0$ if $E(\G)$ does not contain $(i,j).$ We consider such an {\em inverse\/} correspondence in order to model preference digraphs in Section~\ref{leader}: in this case, the transitions in the Markov chain are performed from ``worse'' objects to ``better'' ones, so the Markov chain stochastically ``searches the leaders.''

\begin{theorem}\label{M}
For any finite Markov chain$,$ its matrix of Ces\'aro limiting probabilities coincides with the matrix $\J$ of any digraph inversely corresponding to this Markov chain.
\end{theorem}

Theorem~\ref{M} follows from the {\em Markov chain tree theorem\/} first proved by Wentzell and Freidlin \cite{WentzellFreidlin70a} and rediscovered
in \cite{LeightonRivest83,LeightonRivest86}, which, in turn, can be immediately proved using item~8 of Theorem~\ref{th4} and a result of~\cite{Rothblum76a} (see~\cite{CheAga02a+}). Another proof of Theorem~\ref{M} can be found in \cite{chebotarev02spanning}. A review on forest representations of Markov chain probabilities is given in~\cite{Che04MCTT}. For an interpretation of $J(\tau)$ in terms of Markov chains we refer to~\cite{AgaChe01}. Theorem~\ref{M} will be used in Section~\ref{leader}.

\section{Forest based accessibility measures}
\label{dosti}

Formally, by an {\em accessibility measure\/} for digraph vertices we mean any function that assigns a matrix $P=(p\_{ij})\_{n\times n}$ to every weighted digraph $\G,$ where $n=\card{V(\G)}.$ Entry $p\_{ij}$ is interpreted as the accessibility (or connectivity, relatedness, proximity, etc.) of $j$ from $i.$

Consider the accessibility measures $P^{{\rm out}}_{\tau}=J(\tau),$ where $J(\tau)$ is defined by (\ref{Jtau}), and
$P^{{\rm in}}_{\tau}=(p^{{\rm in}}_{ij})$ with
$p^{{\rm in}}_{ij}=\e(\FF^{i\tor j}(\tau))/\e(\FF^{\tor}(\tau)),$ where
$\FF^{i\tor j}(\tau)$ and
$\FF^{\tor}(\tau)$ are, respectively, the
$\FF^{i\tor j}$ and
$\FF^{\tor}$ for the digraph $\G(\tau)$ obtained from $\G$ by the multiplication of all arc weights by~$\tau.$ Parameter $\tau$ specifies the relative weight of short and long connections\x{} in~$\G$.

\begin{definition}\label{dua}
{\em Accessibility measures $P^{(1)}$ and $P^{(2)}$ are {\em dual\/} if for every $\G$ and every $i,j\in V(\G),\;$ $p^{(1)}_{ij}(\G)=p^{(2)}_{ji}(\G'),$ where $\G'$ is obtained from $\G$ by the reversal of all arcs (preserving their weights).}
\end{definition}

The following proposition results from the fact that the reversal of all arcs in $\G$ transforms all out-forests into in-forests and vice versa.

\begin{prop}\label{pr_dual}
For every $\tau>0,$ the measures $P^{{\rm out}}_{\tau}$ and
                                 $P^{{\rm  in}}_{\tau}$ are dual.
\end{prop}

What is the difference in interpretation between
$P^{{\rm out}}_{\tau}$ and
$P^{{\rm in}}_{\tau}$? A partial answer is as follows.
$P^{{\rm out}}_{\tau}$ can be interpreted as the relative weight of
$i\to j$ connections among the out-connections of $i,$ whereas
$P^{{\rm in}}_{\tau}$ is the relative weight of
$i\to j$ connections among the in-connections of $j.$ Naturally, these relative weights need not coincide. For example, a connection between an average man and a celebrity is usually more significant for the average man. This example demonstrates that self-duality is not an imperative requirement to accessibility measures. Still, self-dual measures are quite common.
Their properties have been studied in~\cite{CheSha98}.

The following conditions some of which were proposed in \cite{CheSha98} can be considered as desirable properties of vertex accessibility measures.

\condition{Nonnegativity}
{$p\_{ij}\ge0,\;\:i,j\in V(\G)$.}

\conditiont{Reachability condition}
{For any $i,j\in V(\G),\;$ ($p\_{ij}=0 \Leftrightarrow j$ is unreachable from~$i$).}

\conditiont{Self-accessibility condition}
{For any distinct $i,j\in V(\G),$ (A)~$p\_{ii}>p\_{ij}$ and (B)~$p\_{ii}>p\_{ji}$ hold.}

\conditiont{Triangle inequalities for proximities} 
{If $i,k,t\in V(\G)$, then
(A)~$p\_{ki}-p\_{ti}$ $\le p\_{kk}-p\_{tk}$ and
(B)~$p\_{ik}-p\_{it}$ $\le p\_{kk}-p\_{kt}.$ 
}

The triangle inequalities for proximities is a counterpart of the ordinary triangle inequality which characterizes distances (cf.~\cite{CheSha98a}).

Let $k,i,t\in V(\G)$. We say that $k$ {\em is a cutpoint between $i$\/ and $t$\/} if\/ $\G$ contains a path from $i$ to $t,$ $i\ne k\ne t,$ and every path from $i$ to $t$ includes~$k.$

\condition{Transit property}
{If $k$ is a cutpoint between $i$\/ and $t,$ then (A)~$p\_{ik}>p\_{it}$ and (B)~$p\_{kt}>p\_{it}.$ }

\conditiont{Monotonicity} {Suppose that the weight $\e\_{kt}$ of some arc $(k,t)$ is increased or a new $(k,t)$ arc is added to $\G$, and $\D p\_{ij},\;i,j\in V(\G),$ are the resulting increments of the accessibility values. Then{\rm:}\\
{\rm (1)} $\D p\_{kt}>0;$\\
{\rm (2)} If $t$ is a cutpoint between $k$\/ and $i$, then $\D p\_{ki}>\D p\_{ti};$ if $k$ is a cutpoint between $i$\/ and $t$ then $\D p\_{it}>\D p\_{ik};$\\
{\rm (3)} (A)~If $t$ is a cutpoint between $k$\/ and $i$, then $\D p\_{kt}>\D p\_{ki};$\\
\phantom{{\rm (3)}} (B)~If $k$ is a cutpoint between $i$\/ and $t$, then $\D p\_{kt}>\D p\_{it}$. }

\conditiont{Convexity} {(A)~If $p\_{ki}>p\_{ti}$ and $i\ne k,$ then there exists a $k$ to $i$ path such that the difference $p\_{kj}-p\_{tj}$ strictly decreases as $j$ advances from $k$ to $i$ along this path. (B)~If $p\_{ik}>p\_{it}$ and $i\ne k,$ then there exists an $i$ to $k$ path such that the difference $p\_{jk}-p\_{jt}$ strictly increases as $j$ advances from $i$ to $k$ along this path.}

The results of testing $P^{{\rm out}}_{\tau}$ and $P^{{\rm in}}_{\tau}$ are collected{\x} in

\begin{theorem}\label{otledostup}
$\!\!\!${\rm~\cite{Cheb0508:Matrices}.}
The measures $P^{{\rm out}}_{\tau}$ and $P^{{\rm in}}_{\tau}$ satisfy
all the above conditions not partitioned into\/ $(A)$
                                          and\/~$(B)$.
Further\-more$,$
$P^{{\rm out}}_{\tau}$ obeys all\/~$(A)$ conditions and
$P^{{\rm in}}_{\tau}$        all\/~$(B)$ conditions.
\end{theorem}

Consider now the accessibility measures
$\widetilde P^{{\rm out}}=(p\_{ij})=\vj
                        =\lim_{\tau\to\infty}P^{{\rm out}}_{\tau}\/$ and
$\widetilde P^{{\rm in}}=\lim_{\tau\to\infty}P^{{\rm  in}}_{\tau}$. Having in mind Theorem~\ref{M}, we call
$\q_{ij}$ the {\it limiting out-accessibility of $j$ from $i$}.

Let us say that a condition {\em is satisfied in the nonstrict form\/} if it is not generally satisfied, but it becomes true after the substitution of $\ge$ for $>,$ $\le$ for $<,$\x{} and ``nonstrictly'' for ``strictly'' in the conclusion of this condition.

Similarly to Proposition~\ref{pr_dual},\x{} we have

\begin{prop}
\label{pr_dual1}
The accessibility measures $\widetilde P^{{\rm out}}$ and
                           $\widetilde P^{{\rm in}}$ are dual.
\end{prop}

The results of testing $\widetilde P^{{\rm out}}$ and $\widetilde P^{{\rm in}}$ are collected{\x} in

\begin{theorem}
\label{bliz}
The accessibility measures $\widetilde P^{{\rm out}}$ and
                           $\widetilde P^{{\rm in}}$
satisfy nonnegativity and the\/ {\rm ``$\Leftarrow$''} part of reachability condition$,$ but they violate the\/ {\rm ``$\Rightarrow$''} part of the reachability condition. Moreover$,$
$\widetilde P^{{\rm out}}$ satisfies$,$ in the nonstrict form$,$ items~$(A)$ of self-accessibility condition$,$ transit property$,$ monotonicity$,$ and convexity$,$ whereas
$\widetilde P^{{\rm in}}$ satisfies$,$ in the nonstrict form$,$ items~$(B)$ of these conditions. The corresponding strict forms are violated.\x{}
$\widetilde P^{{\rm out}}$ satisfies $(A)$ and
$\widetilde P^{{\rm in}}$  satisfies $(B)$ of triangle inequality for proximities.
\end{theorem}

By virtue of Theorem~\ref{bliz}, the limiting accessibility measures only ``marginally'' correspond to the concept of accessibility underlying the above conditions.

\begin{proof}
{The nonstrict satisfaction of the conditions listed in the theorem follows from Theorem~\ref{otledostup}, Proposition~\ref{pr_dual1} and item~5 of Theorem~\ref{th4}. To prove that the strict forms of these conditions and the ``$\Rightarrow$'' part of reachability condition are violated, it suffices to consider the digraph $\G$ with $n\ge3$, $E(\G)=\{(1,2),(2,3)\}$, and $\e\_{12}=\e\_{23}=1$.}
\end{proof}

Let us mention one more class of accessibility measures: those of the form $(I+\a\J)^{-1}$, $0<\a<\si\_{(\di)}/\si\_{(\di+1)}$. These measures are ``intermediate'' between $P^{{\rm out}}_{\tau}$ and $\widetilde P^{{\rm out}}$, because they are positive linear combination of $J\_{(\di)}$ and $J\_{(\di+1)}$~\cite{AgaChe01}. That is why we termed them the {\em matrices of dense out-forests}. In the terminology of \cite[p.~152]{MeyerStadelmaier78}, $(I+\a\J)^{-1}$ with various sufficiently small $\a>0$ make up a class of {\em nonnegative nonsingular commuting weak inverses\/} for~$L$. These measures and the dual measures have been studied in \cite{AgaChe01} (see also \cite[p.~270--271]{CheAga02a+}). Other interesting related topics are the forest distances \cite{chebotarev02forest} and the forest based centrality measures~\cite{CheSha97}.

Thus, we can conclude that forest based accessibility measures may be useful for the analysis of various networks.

\section{Rooted forests and the problem of leaders}
\label{leader}

Ranking from tournaments or irregular pairwise contests is an old, but still intriguing problem. Its statistical version is ranking objects on the basis of paired comparisons~\cite{David88}. Analogous problems of the analysis of individual and collective preferences arise in the contexts of policy, economics, management science, sociology, psychology, etc. Hundreds of methods have been proposed for handling these problems (for a review, see, e.g., \cite{David88,DavidAndrews93,CookKress92,BelkinLevin90,CheSha97a,CheSha99,Laslier97,DiazHendrickxLohmann13SCW,csato2012ranking}).

In this section, we consider a weighted digraph~$\G$ that represents a competition (which need not be a round robin tournament, i.e., can be ``incomplete'') with weighted pairwise results. The digraph can also represent an arbitrary weighted preference relation. The result we present below can be easily extended to multidigraphs.

One of the popular exquisite methods for assigning scores to the participants in a tournament was independently proposed by Daniels \cite{Daniels69}, Moon and Pullman~\cite{MoonPullman69,MoonPullman70}, and Ushakov~\cite{Ushakov71,Ushakov76} (for its axiomatization, see~\cite{SlutzkiVolij05}) and reduces to finding nonzero and nonnegative solutions to the system of equations
\begin{equation}\label{D}
Lx=0.
\end{equation}

Entry $x_i$ of a solution vector $x=(x\_1\cdc x\_n)^{\interca}\x{}$ is considered as a sophisticated ``score'' attached to vertex~$i$. This method was multiply rediscovered with different motivations (some references are given in~\cite{CheSha99}). As Berman \cite{Berman80} noticed (although, in other contexts, similar results had been obtained by Maxwell \cite{Maxwell1892} and other writers, see~\cite{CaplanZeilberger82}), if a digraph is strong, then the general solution to (\ref{D}) is provided by the vectors proportional to $t=(t\_1\cdc t\_n)^{\interca},$ where $t\_j$ is the weight of the set of spanning trees (out-arborescences) diverging from~$j$. This fact can be easily proved as follows. By the matrix-tree theorem for digraphs (see, e.g.,~\cite{Harary69}), $t\_j$ is the cofactor of any entry in the $j$th column of~$L$. Then for every $i\in V(\G),$ $\sum^n_{j=1}\,\l\_{ij}\,t_j=\det L$ (the row expansion of $\det L$) and, since $\det L=0,\;$ $t$ is a solution to~(\ref{D}). As $\rank L=n-1$ (since the cofactors of $L$ are nonzero), any solution to (\ref{D}) is proportional to~$t$.

Berman \cite{Berman80} and Berman and Liu~\cite{BermanLiu96} asserted that this result is sufficient to rank the players in an arbitrary competition, since the strong components of the corresponding digraph supposedly ``can be ranked such that every player in a component of higher rank defeats every player in a component of lower rank. Now by ranking the players in each component we obtain a ranking of all the players.'' While the statement about the existence of a natural order of the strong components is correct in the case of round-robin tournaments, it need not be true for arbitrary digraphs that may have, for instance, several source knots. That is why, the solution devised for strong digraphs does not enable one to rank the vertices of an arbitrary digraph.

Let us consider the problem of interpreting, in terms of forests, the general solution to (\ref{D}) and the problem of choosing a particular solution that could serve as a reasonable score vector in the case of arbitrary digraph~$\G$.

If $\G$ contains more than one source knot, there is no spanning diverging tree in~$\G$. Recall that $K\_1\cdc K\_{\di}$ are the source knots of $\G,$ where $\di$ is the out-forest dimension of $\G,$ and $\ktil=\bigcup_{i=1}^{d'}K\_i$.

Suppose, without loss of generality, that the vertices of $\G$ are numbered as follows. The smallest numbers are attached to the vertices in $K\_1$, the following numbers to the vertices in $K\_2$, etc., and the largest numbers to the vertices in $V(\G)\setminus\ktil.$ Such a numeration we call {\em standard}.

\begin{theorem}\label{d_sol}
Any column of $\J$ is a solution to~$(\ref{D})$.
Suppose that the numeration of vertices is standard and $j\_1\in K\_1\cdc j\_{\di}\in K\_{\di}$. Then the columns $\J\_{\bullet j\_1}\cdc\J\_{\bullet j\_{\di}}$ of\/ $\G$ make up an orthogonal basis in the space of solutions to~$(\ref{D})$ and
$\J\_{\bullet j\_s}
=\e^{-1}(\TT\_s)\bigl(0\cdc 0,\e(\TT_s^{i\_s+1})\cdc$ $\e(\TT_s^{i\_s+k\_s}),0\cdc 0\bigr)^{\interca},$
where $\{i\_s+1\cdc i\_s+k\_s\}=K\_s$ and $\TT_s$ is the set of out-arborescences of $K\_s,$ $s=1\cdc \di.$
\end{theorem}

By virtue of Theorem~\ref{d_sol}, the general solution to (\ref{D}) is the set of all linear combinations of partial solutions that correspond to each source knot of~$\G$.

\begin{proof}
{The first statement follows from $L\q=0$ (item~3 of Theorem~\ref{th4}). By item~6 of Theorem~\ref{th4}, $\rank\q=\di$ and $\rank L=n-\di$. Hence, $\di$ is the dimension of the space of solutions to~(\ref{D}). Let $j\_s\in K\_s,\;s=1\cdc\di.$ Then, by items~1 and~2 of Theorem~\ref{th2},
$$
\J\_{\bullet j\_s}
=\e^{-1}(\TT\_s)
\bigl(0\cdc 0,\e(\TT_s^{i\_s+1})\cdc\e(\TT_s^{i\_s+k\_s}),
      0\cdc 0\bigr)^{\interca}.
$$
These $\di$ solutions to~(\ref{D}) are orthogonal and thus, linearly independent.
}
\end{proof}

As a reasonable ultimate score vector, the arithmetic mean $x={1\over n}\J\!\cdot\!(1\cdc 1)^{\interca}$ of the columns of $\q$ can be considered. A nice interpretation of this vector is given by

\begin{corollary}\label{c-lim_distr}
{\rm (of Theorem~\ref{M}).}
For any Markov chain inversely corresponding to $\G,$ $x={1\over n}\J\!\cdot\!(1\cdc 1)^{\interca}$ is the limiting state distribution$,$ provided that the initial state distribution is uniform.
\end{corollary}

It can be mentioned, however, that the ranking  method based on $\J$ takes into account long paths in $\G$ only. That is why, in any solutions to (\ref{D}), the vertices that are not in the source knots are assigned zero scores, which is questionable. The estimates based on the matrices $Q(\tau),$ instead of $\q,$ are free of this feature. On the other hand, both methods violate the {\em self-consistent monotonicity\/} axiom~\cite{CheSha99}, and so do the methods that count the {\it walks\/} between vertices. This axiom is satisfied by the {\it generalized Borda method\/} \cite{Che89,Che94} that produces the score vectors $J'(\tau)\!\cdot\!(\od(1)-\id(1)\cdc\od(n)-\id(n))^{\interca}$, where $J'(\tau)$ is the matrix $J(\tau)$ of the undirected graph corresponding to $\G$~\cite{Sha94}. In our opinion, the latter method can be recommended as a well-grounded approach to scoring objects on the base of arbitrary weighted preference relations, incomplete tournaments, irregular pairwise contests, etc.

\section{A concluding remark: a communicatory interpretation of some forest matrices}

In closing, let us mention an interpretation of forest matrices in terms of information dissemination. Consider the following metaphorical model. First, a plan of information transmission along a digraph is chosen. Such a plan is a diverging forest $F\in\FF^{\rto}$: the information is injected into the roots of $F$; then it ought to come to the other vertices along the arcs of $F$. Suppose that $\e\_{ij}\!\in{]}\xy 0,1]$ is the probability of successful information transmission along the $(i,j)$ arc, $i,j\in V(\G),$ and that the transmission processes in different arcs are statistically independent. Then $\e(F)$ is the probability that plan $F$ is successfully realized. Suppose now that each plan is selected with the same probability $\card{\FF^{\rto}}^{-1}.$ Then $J\_{ij}$ (see~(\ref{J})) is the probability that the information came to $j$ from root $i$, provided that the transmission was successful. As a result, if one knows that the information was corrupted at root $i$ and the transmission was successful, then $J\_{ij}$ is the probability that this corrupted information came to~$j.$

Similarly, interpretations of this kind can be given to other stochastic forest matrices. This model is compatible with that of centered partitions \cite{Lenart98} and comparable with certain models of~\cite{Pavlov00}.

\bibliographystyle{elsarticle-num}
\bibliography{all2}
\end{document}